\documentclass[review,hidelinks,onefignum,onetabnum]{siamart220329}

\usepackage{amsfonts,bigints,graphicx,caption,subcaption,xcolor}

\newcommand{\C}{\mathbb{C}}
\newcommand{\R}{\mathbb{R}}
\newcommand{\Et}{\boldsymbol{\eta}}
\newcommand{\W}{\boldsymbol{\wp}}
\newcommand{\Ze}{\boldsymbol{\zeta}}
\newcommand{\Z}{\boldsymbol{Z}}
\newcommand{\z}{\boldsymbol{z}}
\newcommand{\hz}{\hat{z}}
\newcommand{\sz}{\mathsf{z}}
\renewcommand{\a}{\mathfrak{a}}
\renewcommand{\b}{\mathfrak{b}}
\renewcommand{\c}{\mathfrak{c}}
\renewcommand{\d}{\mathfrak{d}}
\newcommand{\zz}{\mathfrak{z}}
\newcommand{\CC}{\mathfrak{C}}
\newcommand{\GG}{\mathfrak{G}}
\newcommand{\ZZ}{\mathfrak{Z}}
\newcommand{\mR}{\mathcal{R}}
\newcommand{\af}{\alpha}
\newcommand{\ta}{\tan\af}
\newcommand{\ca}{\cot\af}
\newcommand{\eps}{\epsilon}

\newcommand{\BE}{\begin{equation}}
\newcommand{\EE}{\end{equation}}
\renewcommand{\o}{\omega}
\newcommand{\deh}{\partial}
\newcommand{\ov}{\overline}

\newsiamremark{remark}{Remark}
\newsiamremark{hypothesis}{Hypothesis}
\crefname{hypothesis}{Hypothesis}{Hypotheses}
\newsiamthm{claim}{Claim}

\headers{The Green's Function on Rhombic Flat Tori}{A. E. D. Castillo, G. A. Lobos, and V. Ramos Batista}

\title{The Green's Function on Rhombic Flat Tori\thanks{
\funding{The first author is supported by the Coordination of Superior Level Staff Improvement (CAPES)} through his master's scholarship No. 88887.976809/2024-00.}}

\author{A. E. D. Castillo\thanks{Department of Mathematics, Federal University of S\~ao Carlos, S\~ao Carlos - SP, Brazil (\email{axel.diaz@estudante.ufscar.br}, \email{lobos@ufscar.br}, \url{https://www.dm.ufscar.br/profs/lobos}).}
\and G. A. Lobos\footnotemark[2]
\and V. Ramos Batista\thanks{Centre of Mathematics, Cognition and Computer Science, Federal University of ABC, Santo Andr\'e - SP, Brazil (\email{valerio.batista@ufabc.edu.br}, \url{https://www.ufabc.edu.br/ensino/docentes/valerio-ramos-batista}).}}

\usepackage{amsopn}

\begin{document}
\nolinenumbers
\maketitle

\begin{abstract}
We obtain the Green's function $G$ for any flat rhombic torus $T$, always with numerical values of significant digits up to the fourth decimal place (noting that $G$ is unique for $|T|=1$ and $\int_TGdA=0$). This precision is guaranteed by the strategies we adopt, which include theorems such as the Legendre Relation, properties of the Weierstra\ss\,P-Function, and also the algorithmic control of numerical errors. Our code uses complex integration routines developed by H. Karcher, who also introduced the symmetric P-Weierstra\ss\,function, and these resources simplify the computation of elliptic functions considerably.
\end{abstract}



\section{Introduction}
The Green's function $G$ is one of the single most powerful resources to solve Partial Differential Equations, specially those of Electromagnetism, Acoustics, Hydrodynamics and Particle Physics, among others. Although most of its applications occur in Physics, this function is an important mathematical tool, and in many cases it can be studied exclusively with methods in Geometry and Complex Analysis.

\begin{figure}[htbp]
\center
\includegraphics[scale=0.31]{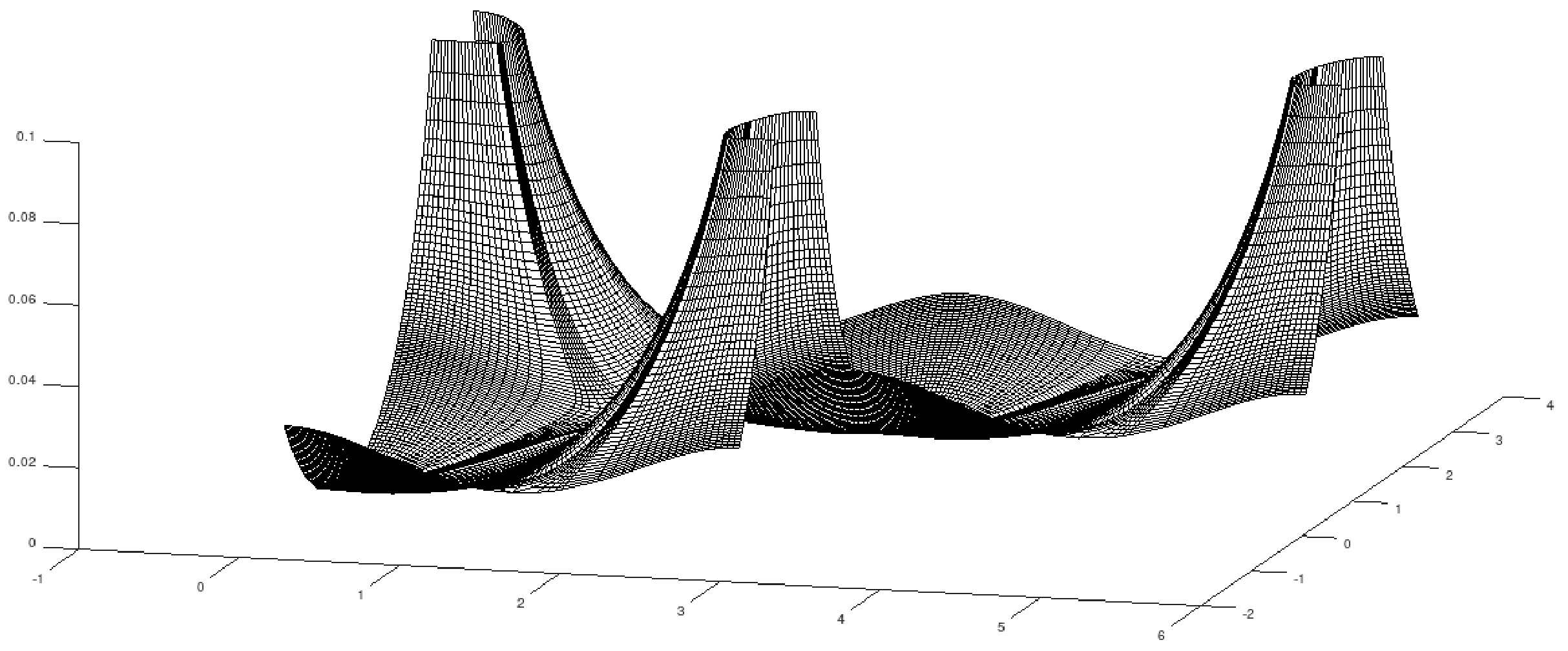}
\smallskip
\caption{A piece of the graph of $G$ for $\rho=\pi/4$ (see later).}
\label{fig:graph}
\end{figure}

By starting from the theoretical definition of $G$, which is a second-order boundary-value problem on a certain domain, the first difficulty is to solve it by an explicit formula. Since this paper focus on surface domains, namely flat rhombic tori $T$, then $G$ is the unique solution of
\BE
   \Delta G = -\delta+\frac{1}{|T|}, \int_TGdA=0,\label{eq:def_G}
\EE
where $\Delta$ is the Laplacian on $T$ and $\delta$ the Dirac delta distribution at $0\in\C$, if we take the complex plane as the universal covering of $T$. Of course, on $\C$ both $G$ and $\delta$ turn out to be doubly periodic, and our numerical simulations always adopt this approach, as exemplified in \cref{fig:graph}.

In point of fact, this figure depicts what we call {\it the non-negative $G$}. Namely, the one that vanishes at its global minimum, and so $\int_TGdA$ is a positive constant. However, for many practical purposes one does not really need this constant. This is the case, for instance, of any linear differential operator whose homogeneous equation admits only the trivial solution. Nevertheless, we shall present a way to compute $\int_TGdA$ later on.

Other surface domains, like the ellipse and the ellipsoid, where fully investigated in \cite{frank1,frank2} for the purpose of studying the equilibrium configuration of self-inhibitory droplets trapped there. Recent works like \cite{L,R,S} have already applied $G$ on flat tori, but there the authors only looked for function minima without obtaining a fully numerical $G$.

Once having obtained $G$ through an explicit formula, the second difficulty is to compute numerical values with accuracy. The reader might even consider it as a minor problem, since present-day computers work with double precision. For instance, GNU-Octave represents any floating-point number with 64-bits (in the IEEE 754 format), hence a {\it few} number of arithmetic operations will give results still accurate at least in single precision, namely 32-bits. By the way, in this paper we use 3.6GB of RAM, microprocessor Intel(R) Pentium(R) Gold 7505 \@ 2.00GHz, operating system Linux Ubuntu 20.04.6, and Octave 5.2.0 (R2020). With this setup we invoke {\tt format long} and compute {\tt 10\^{}5*pi-314159} to get {\tt 2.65358979\underline{2898856}e-01}, underlined to show that $\pi$ is accurate up to the 14th decimal. The remaining decimals are due to {\it roundoff error}, because even quite simple fractions in base 10 turn out to be repeating decimals in base 2. The easiest example is $1/5$, represented as $0.\ov{0011}$ and rounded up by the computer to a slightly greater number, approximately $0.2+10^{-17}$. Another example is the $\eps$-machine. In Octave we get {\tt eps} as {\tt 2.2204e-16}, so that the expression {\tt 1+eps>1} is true, but {\tt 1+eps/2>1} is false. Namely, many terms and factors will be neglected by others whose order of magnitude prevail. 

We have just showed examples with at most three arithmetic operations, but $G$ can be given either by an integral of elliptic functions, or an infinite sum of terms or product of factors. Of course, the computer will deal with a finite approximation, which in any case will lead to {\it method}, {\it truncation} and {\it roundoff} errors. There are various references to them, but here we cite \cite{BF}. In order to get reliable results, one could resort to computationally expensive routines, but the machine error will accumulate and propagate for exceedingly many arithmetic operations. Computational errors used to be tracked more cautiously before 1998, when the quadruple precision IEEE 754 started. However, personal computers still adopt the double precision simply because of the law of supply and demand. This means, researchers who really need more accuracy will resort to supercomputers instead of desktops, and they can easily adapt our methodology for their purposes.

As a matter of fact, controlling computational errors used to be feasible until the early 90s, and for the readers who want to have a bird's-eye view we cite \cite[Chapter~10]{D}. However, with time one had to resort to other strategies, since numerical computations started involving more elaborate algorithms. Here we cite H. Karcher's complex integration routines, available at \url{https://www.math.uni-bonn.de/people/karcher}, in which the author replaced the classical Simpson's integration method with 4th degree polynomials. At that time, Karcher produced many 2-D views of minimal surfaces in \cite{K1} with BASIC programming language. But he kept control of important surface data, like the mean and Gaussian curvatures, which can be computed for discrete surfaces (see \cite[Chapter~16]{B} for details). At that time, one avoided overloading the program with heavy computations, and so he found out that Simpson's approach, or even cubic splines, were insufficient to get reliable results. By testing the numerical integration with 4th degree polynomials, he finally came out with the sought after minimal surfaces.

In this paper we follow a similar approach in order to obtain the Green's function $G$ for any flat rhombic torus $T$, always with numerical values of significant digits up to the fourth decimal place (noting that $G$ is unique for $|T|=1$ and $\int_TGdA=0$). This precision is guaranteed by the strategies we adopt, which include theorems such as the Legendre Relation, properties of the Weierstra\ss\,P-Function, and also the algorithmic control of numerical errors. Our code uses complex integration routines developed by Karcher, who also introduced the symmetric P-Weierstra\ss\,function, and these resources simplify the computation of elliptic functions considerably.

The paper is organised as follows. Our main results are in \cref{sec:main}, our programming strategies described in \cref{sec:alg}, experimental results in \cref{sec:experiments}, and the conclusions follow in \cref{sec:conclusions}.

\section{Main results}
\label{sec:main}

We begin with some preliminaries. Let $z$ be the standard variable of the complex plane $\C$, and consider two elements $\o_1$, $\o_2$ in $\C$ that are linearly independent over $\R$. They generate the {\it period lattice} $\Lambda:=\{2m\o_1+2n\o_2:m,n\in\Z\}$. Henceforth we shall indicate any set minus the origin by an asterisk, thus $\Lambda^*=\Lambda\setminus\{0\}$. The {\it classical Weierstra\ss\,P-function} is defined as
\BE
   \W(z)=\frac{1}{z^2}+\sum_{\lambda\in\Lambda^*}\biggl(\frac{1}{(z-\lambda)^2}-\frac{1}{\lambda^2}\biggl),\label{eq:classical_WP}
\EE
which converges uniformly and absolutely on any compact subset of $\C\setminus\Lambda$.

\begin{figure}[htbp]
\center
\includegraphics[scale=0.64]{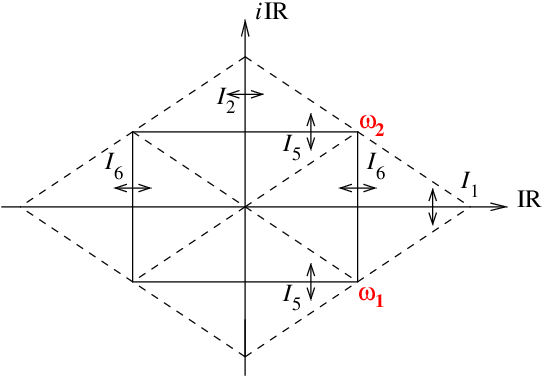}\hfill
\includegraphics[scale=0.72]{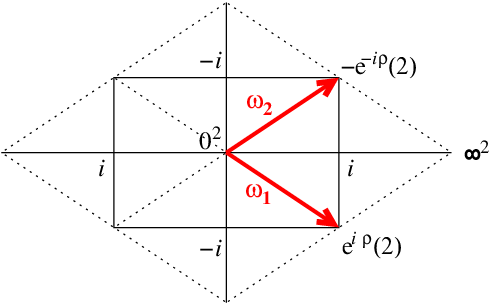}
\centerline{(a)\hspace{6.5cm}(b)}
\smallskip
\caption{(a) Some anti-holomorphic involutions of $T$; (b) some values of $\wp$ on $T$.}
\label{fig:scheme}
\end{figure}

Let $\hz:=z-\o_1-\o_2$, hence $\W(\hz)$ takes the values $\infty$, $e_1$, $e_2$, $e_3$ when $z$ is $\o_1+\o_2$, $\o_2$, $\o_1$ and $0$, respectively. Consider $\c:=\sqrt{(e_1-e_3)(e_3-e_2)}$ such that $\c=-i|\c|$. Then {\it Karcher's symmetric Weierstra\ss\,P-function} is defined as
\BE
   \wp(z)=\frac{\W(\hz)-e_3}{\c},\label{eq:rel_wps}
\EE
for all $z\in\C$. Though quite simple in aspect, Equation \cref{eq:rel_wps} was Karcher's greatest key to construct many examples of minimal surfaces in \cite{K1,K2,K3,K4}, among others. His works boosted the theory of minimal surfaces, and were cited by many authors after him. Now we show the noticeable simplifications of working with $\wp$. \cref{fig:scheme}(a) depicts a typical rhombus with four anti-holomorphic involutions, so that $I_3=I_1\circ I_2$ and $I_4=I_1\circ I_1$ are $180^\circ$-rotation about the origin and the identity, respectively.

For $\rho\in(-\pi/2,\pi/2)$ we follow the same approach of \cite{K1} and represent strategic values of $\wp$ directly on $T$. In \cref{fig:scheme}(b) one sees that $\wp(\o_1)=e^{i\rho}$ with branch order 2, $\wp(\frac{\o_1+\o_2}{2})=i$, $\wp(0)$ is a double zero, and so on. Namely, through symmetries $\wp$ is completely determined by its restriction to one eighth of $T$. \cref{tab:foo} summarises this geometrical simplification, which will be endowed by the primitive $\zeta(z):=-\int_0^z\wp$, as detailed in \cref{sec:alg}.

\begin{table}[htbp]
\footnotesize
\caption{$\wp$ and involutions of $T$.}\label{tab:foo}
\begin{center}
\begin{tabular}{|c|c|c|c|}\hline
Composition & $\wp\circ I_1=\wp\circ I_2$ & $\wp\circ I_3=\wp\circ I_4$ & $\wp\circ I_5=\wp\circ I_6$\\ \hline
Result      & $-\bar{\wp}$               & $\wp$                       & $1/\bar{\wp}$ \\ \hline
\end{tabular}
\end{center}
\end{table}

\

Define $\hat{\R}:=\R\cup\{\infty\}$, $\a=\Re\,\o_1$, $\b=\Im\,\o_2$, $x=\Re\,z$ and $y=\Im\,z$. We obtain:

\begin{theorem} Let $\ZZ:\C\to\hat{\R}$ given by
\BE
   \ZZ(\sz)=-\frac{1}{2\pi}\Re\int_0^\sz\c\zeta(z)dz.\label{def_ZZ}
\EE
Then there exist unique constants $A$, $B$, $C$ depending only on $(\a,\b)$ such that $\hat{\ZZ}:\C\to\hat{\R}$ defined by $\hat{\ZZ}(z):=\ZZ(z)+Ax^2+By^2+Cxy$ is a real-analytic doubly-periodic function.
\label{thm:mainthm}
\end{theorem}

As we are going to see in \cref{sec:alg}, this implies

\begin{corollary}\label{cor:a}
For $A$ and $B$ given in \cref{thm:mainthm} we have $4(A+B)=1/|T|$.
\end{corollary}

Finally, \cref{cor:a} will result in

\begin{theorem}
Consider $\ZZ$ as in \cref{thm:mainthm} and $G=2\ZZ+Ax^2+By^2$. Then $G$ is the Green's function with $G(0)=0$ and $\Delta G=-\sum_{m,n}\delta(2\a+2m\o_1+2n\o_2)+1/|T|$ in $\C$.
\label{thm:def_G}
\end{theorem}

We remark that $\Delta(2\ZZ)$ does {\it not} consist solely of $\delta$-functions as \cref{cor:a} might suggest. The factor 2 will be explained at the end of \cref{sec:experiments}.

\section{Programming strategies}
\label{sec:alg}

\begin{figure}[htbp]
\center
\includegraphics[scale=0.4]{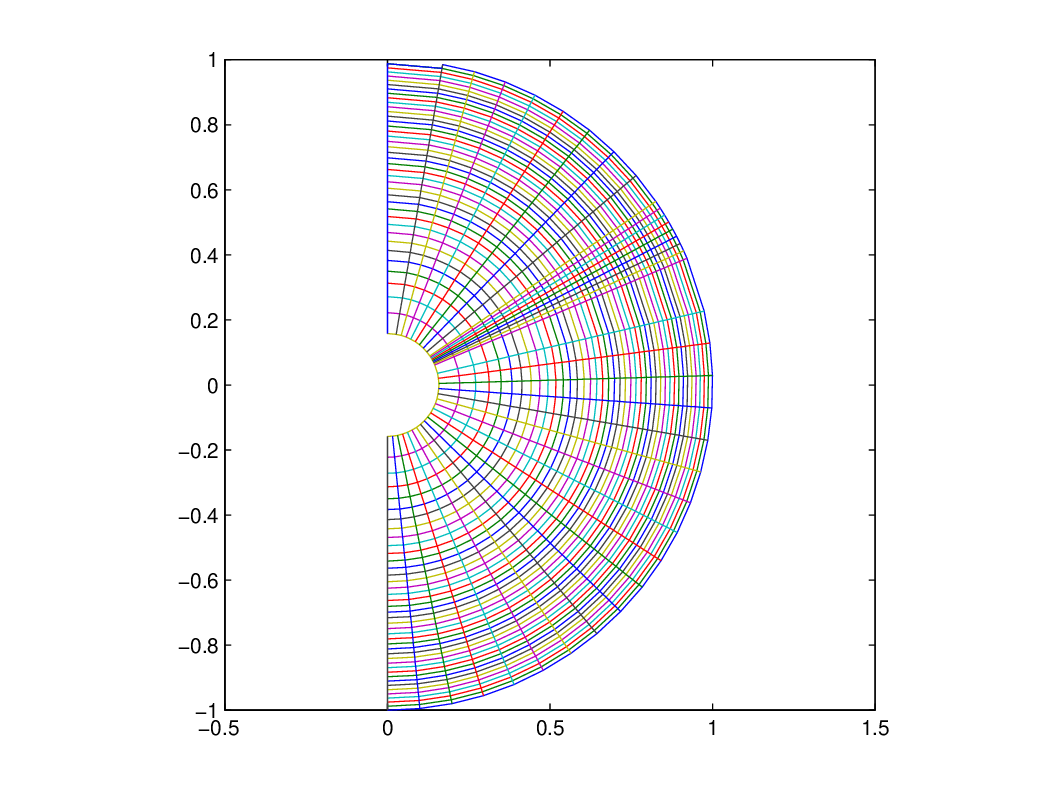}
\includegraphics[scale=0.3]{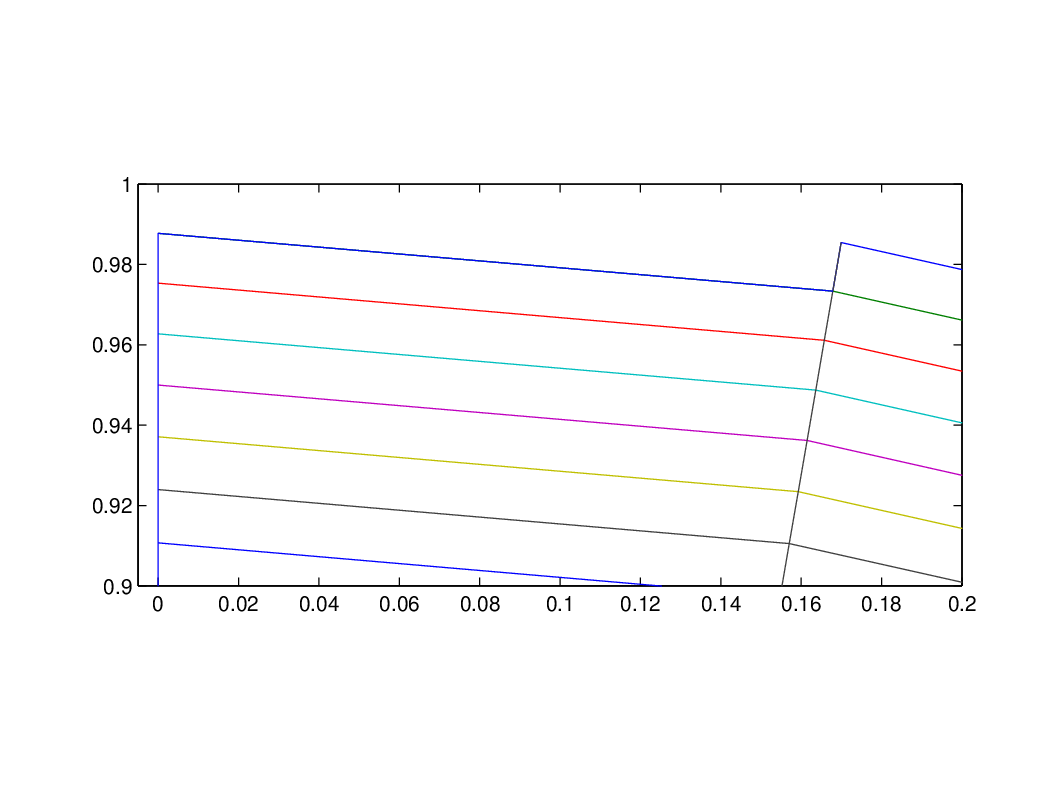}
\centerline{(a)\hspace{6.5cm}(b)}
\smallskip
\caption{(a) $\wp$ for the sub-rectangle of \cref{fig:scheme} with diagonal $\o_1$, (b) detour at $i$.}
\label{fig:WP_img}
\end{figure}

In what follows we shall closely use~\cite[pp.227-232]{A}. \cref{fig:WP_img}(a) represents the image of $\wp$ applied to the sub-rectangle of \cref{fig:scheme} with diagonal $\o_1$. In this example $\rho=1/2$. Notice that we detoured from the origin because our program makes use of $1/\wp$. Another detour, indicated in \cref{fig:WP_img}(b), will soon be made evident by Equation \cref{eq:1st_transf}.

\begin{figure}[htbp]
\center
\includegraphics[scale=0.28]{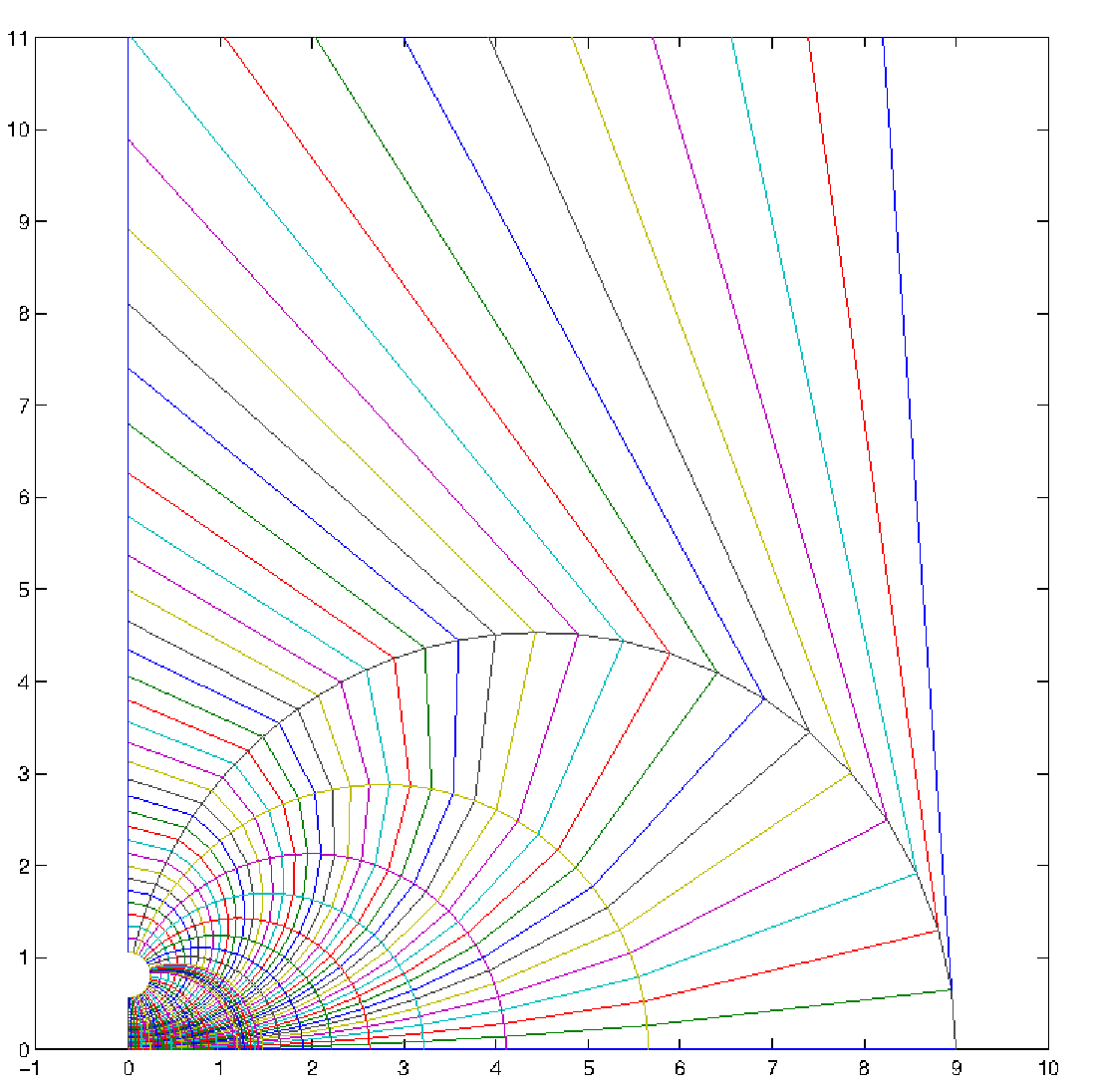}\hfill
\includegraphics[scale=0.43]{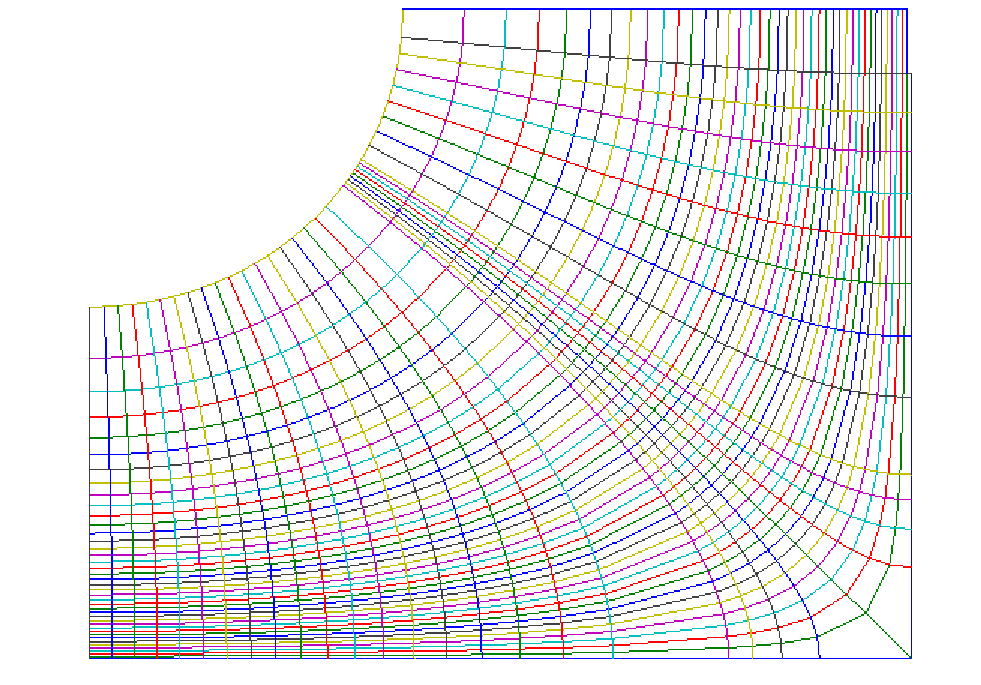}
\centerline{(a)\hspace{6.5cm}(b)}
\caption{(a) Image of $\z(\zz)$, (b) its composition with Equation~\cref{eq:go_torus}.}
\label{fig:fund_rectg}
\end{figure}

Let us call $\zz$ the complex variable in $\C$ restricted to the half-disk of \cref{fig:WP_img}(a). For $\af:=\pi/4+\rho/2$ the M\"obius transformation
\BE
   \z=\sqrt{\ca}\cdot\frac{1-i\zz}{\zz-i}\label{eq:1st_transf}
\EE
takes the half-disk of \cref{fig:WP_img}(a) to the 1st quadrant of $\C$, as depicted in \cref{fig:fund_rectg}(a). We see that $\z^2$ will give the upper-half complex plane, and so Eq.(5) of~\cite[p.230]{A} is re-written as
\BE
   \mathfrak{F}(w)=2\int_0^w\frac{d\z}{[(\ca+\z^2)(\ta-\z^2)]^{1/2}}.\label{eq:go_torus}
\EE

However, Karcher's numerical integration does {\it not} always start at $\z=0$. In general we have
\[
   \mathfrak{F}(w)=2\int_{\z_0}^w\frac{d\z}{[(\ca+\z^2)(\ta-\z^2)]^{1/2}}+{\rm constant},
\]
where the constant is given by
\[
   2\int_0^{\z_0}\frac{d\z}{[(\ca+\z^2)(\ta-\z^2)]^{1/2}}.
\]

Notice that we obtained Equation \cref{eq:go_torus} by applying $\z^2$ to Eq.(5) of~\cite[p.230]{A}, which removed the integrand singularity at the origin. This is crucial for numerical integration methods, because they always rely on {\it uniformly continuous} integrands. In order to control the {\it method error} we first ensured that the integrand of Equation~\cref{eq:go_torus} does not have a singularity at $\z=0$. We can get rid of another two singularities by replacing $\z$ in Equation~\cref{eq:go_torus} with
\BE
   \z=\sqrt{\ta}-\biggl(\sqrt{\sqrt{\ta}-i\sqrt{\ca}}-\Z^2\biggl)^2,\label{eq:get_rid}
\EE
whence
\BE
   \mathfrak{F}(w)=8\bigintsss_{w_0}^{w_1}
   \frac{d\Z}{\big[(\z(\Z)+i\sqrt{\ca})(\z(\Z)+\sqrt{\ta})
   (2\sqrt{\sqrt{\ta}-i\sqrt{\ca}}-\Z^2)\big]^{1/2}}.\label{eq:nice_int}
\EE

The extremes $w_0$ and $w_1$ are left implicit because we shall not use them directly in the numerical integrations.

\begin{figure}[ht!]
\centering
\begin{minipage}{.5\textwidth}
\includegraphics[scale=0.27]{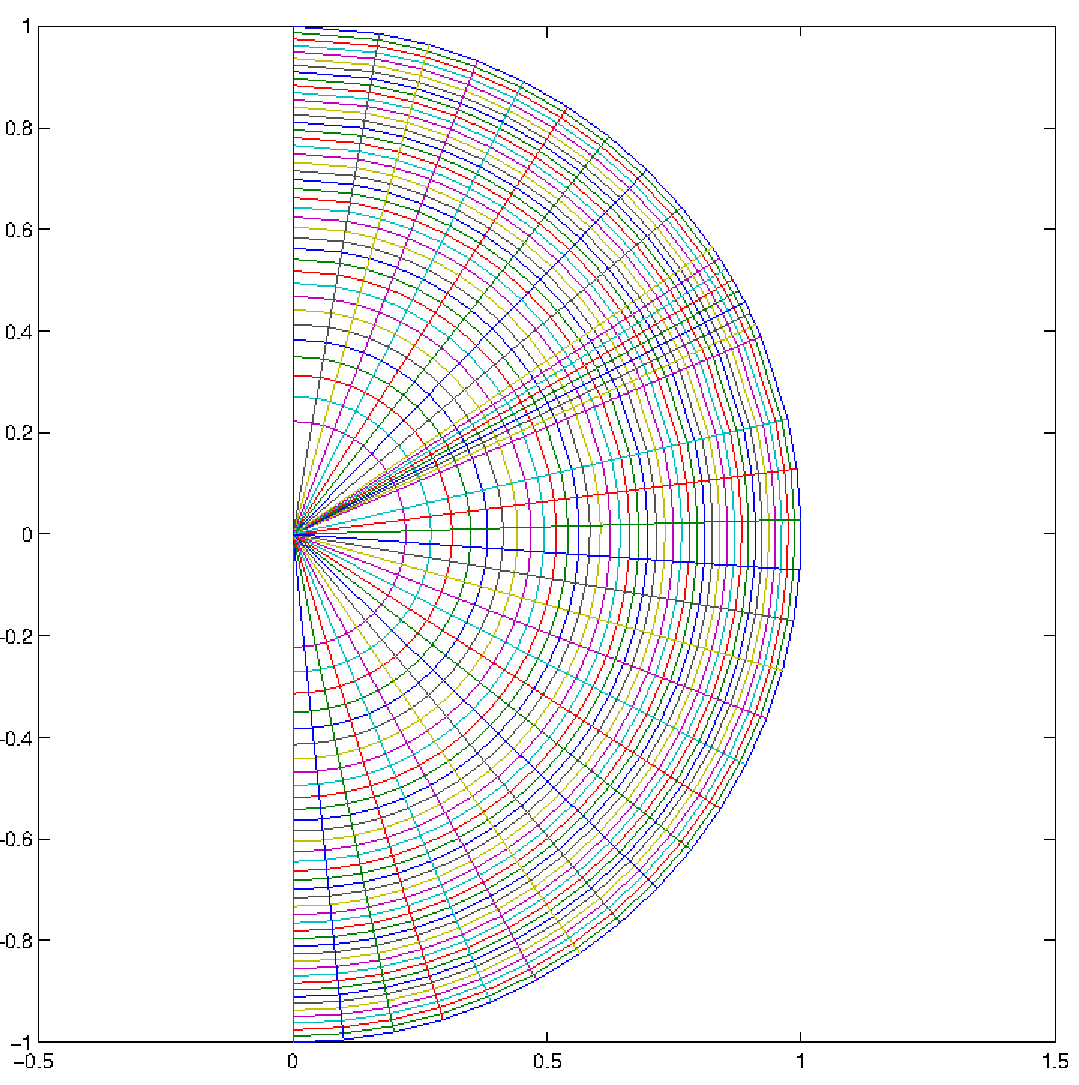}
\caption{Correction of \cref{fig:WP_img}(a).}
\label{fig:fund_rectg_cr}
\end{minipage}%
\begin{minipage}{.5\textwidth}
\includegraphics[scale=0.27]{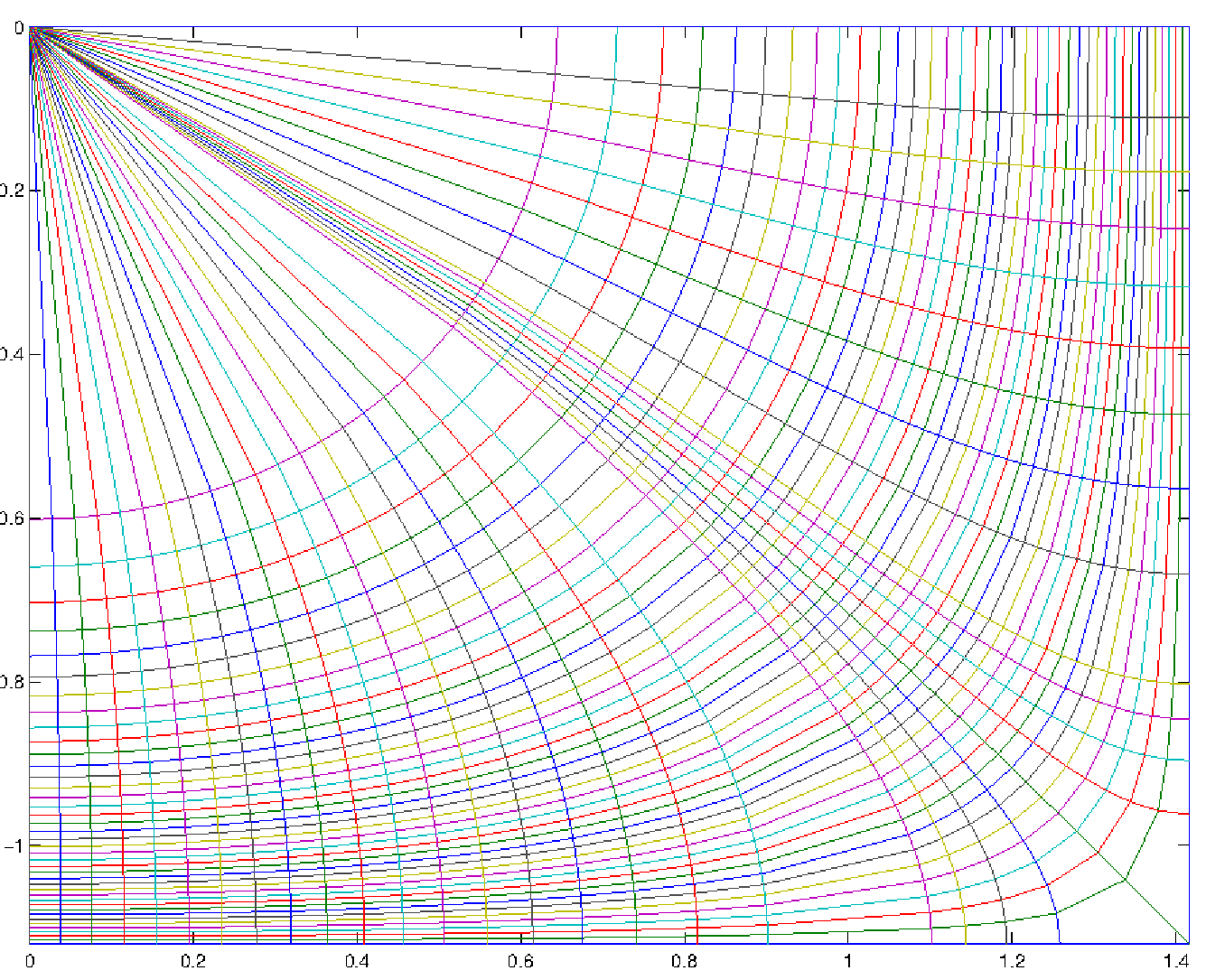}
\caption{Its image under $\wp^{-1}$.}
\label{fig:WP_img_cr}
\end{minipage}
\end{figure}

\cref{fig:fund_rectg}(b) shows Equation \cref{eq:nice_int} applied over the grid in \cref{fig:fund_rectg}(a). Notice that we get a rectangle with a slight detour at the top left corner where $\z=\infty$, which of course cannot belong to \cref{fig:fund_rectg}(a). But the integrand of Equation \cref{eq:nice_int} is free of singularities, even at $\z=\infty$, whence we correct it at the origin. The final result is depicted in \cref{fig:fund_rectg_cr} and \cref{fig:WP_img_cr}.

At this point it is important to notice that 
\BE
   \a=2\int_0^{\sqrt{\ta}}\!\!\!\!\!\frac{d\z}{[(\ca+\z^2)(\ta-\z^2)]^{1/2}},\,\,\,
   \b=2\int_0^{i\sqrt{\ca}}\!\!\!\!\!\frac{-id\z}{[(\ca+\z^2)(\ta-\z^2)]^{1/2}}\label{eq:a_and_b}
\EE
will give the fundamental lattice periods $\o_1:=\a-i\b$ and $\o_2:=\a+i\b$. Indeed, the key difference between the rectangular and rhombic tori is that $\zz=\pm i$ are {\it not} branch points. Hence $\o_{1,2}$ represent {\it diagonals} of a rectangle instead of {\it sides}. Another important remark is that in both Equations \cref{eq:get_rid} and \cref{eq:nice_int} we have $\sqrt{\sqrt{\ta}-i\sqrt{\ca}}$, a term that is computed in our code as {\tt-sqrt(st-i*sc)}, where {\tt st} and {\tt sc} stand for $\sqrt{\ta}$ and $\sqrt{\ca}$, respectively. By the way, the source-codes will be released in future and the reader will be able to access these and other details in full, which for now are only discussed throughout this paper.

From the general theory of elliptic functions we know that $e_1$, $e_2$, $e_3$ are the roots of $4t^3-g_2t-g_3=0$ such that
\[
   g_2=\frac{15}{4}\sum_{m,n}(m\o_1+n\o_2)^{-4}
   \hspace{1cm}{\rm and}\hspace{1cm}
   g_3=\frac{35}{16}\sum_{m,n}(m\o_1+n\o_2)^{-6},
\]
where summations will be always for all $(m,n)\in\mathbb{Z}\times\mathbb{Z}\setminus\{(0,0)\}$ throughout this text. The classical $\Ze$-function is defined as $\Ze(\hz)=-\int^{\hz}\W$ with $\Ze(\hz)-1/\hz\to 0$ when $\hz\to0$. We are going to take $\zeta(z):=-\int_0^z\wp$ as the {\it symmetric zeta-function}. If we consider $\hz\mapsto z$ as a {\it change of variables} then
\[
   \Ze(z)=-e_3z-\c\cdot\int^z\wp
   \hspace{1cm}{\rm s.t.}\hspace{1cm}
   \lim_{z\to\o_1+\o_2}\biggl(\Ze(z)-\frac{1}{z-\o_1-\o_2}\biggl)=0.
\]

\begin{figure}[ht!]
\centering
\begin{minipage}{.5\textwidth}
\includegraphics[scale=0.35]{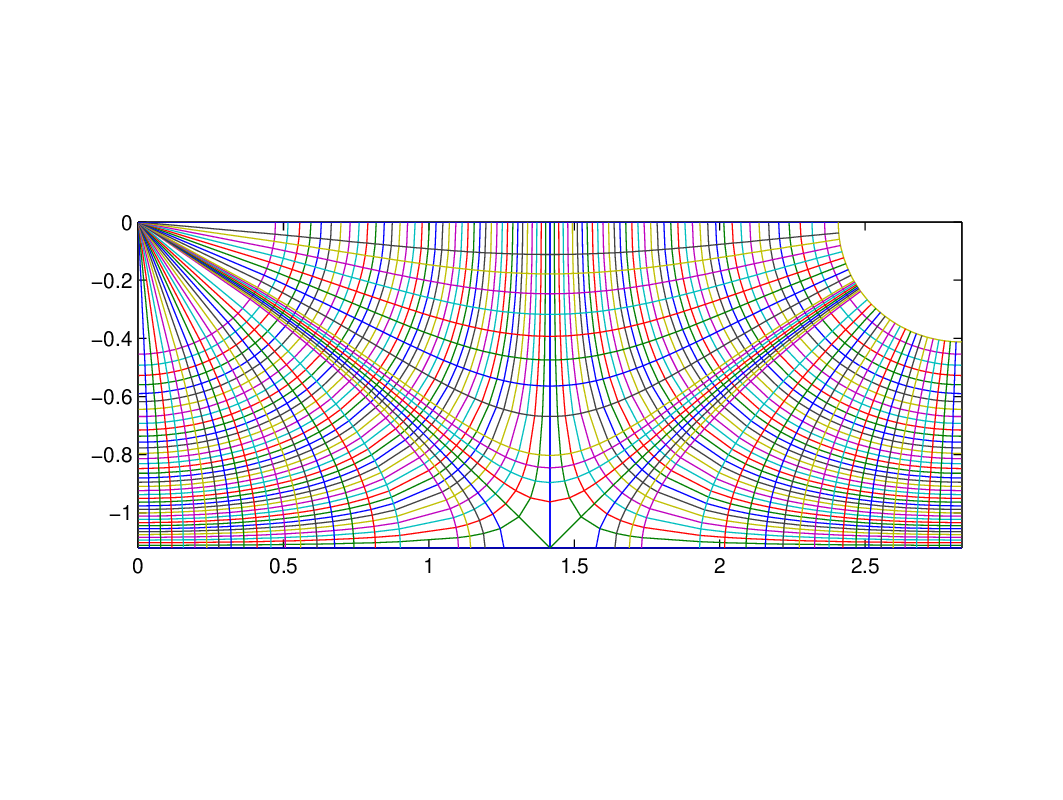}
\caption{Two mirror copies of \cref{fig:WP_img_cr}.}
\label{fig:quarter}
\end{minipage}%
\begin{minipage}{.5\textwidth}
\includegraphics[scale=0.35]{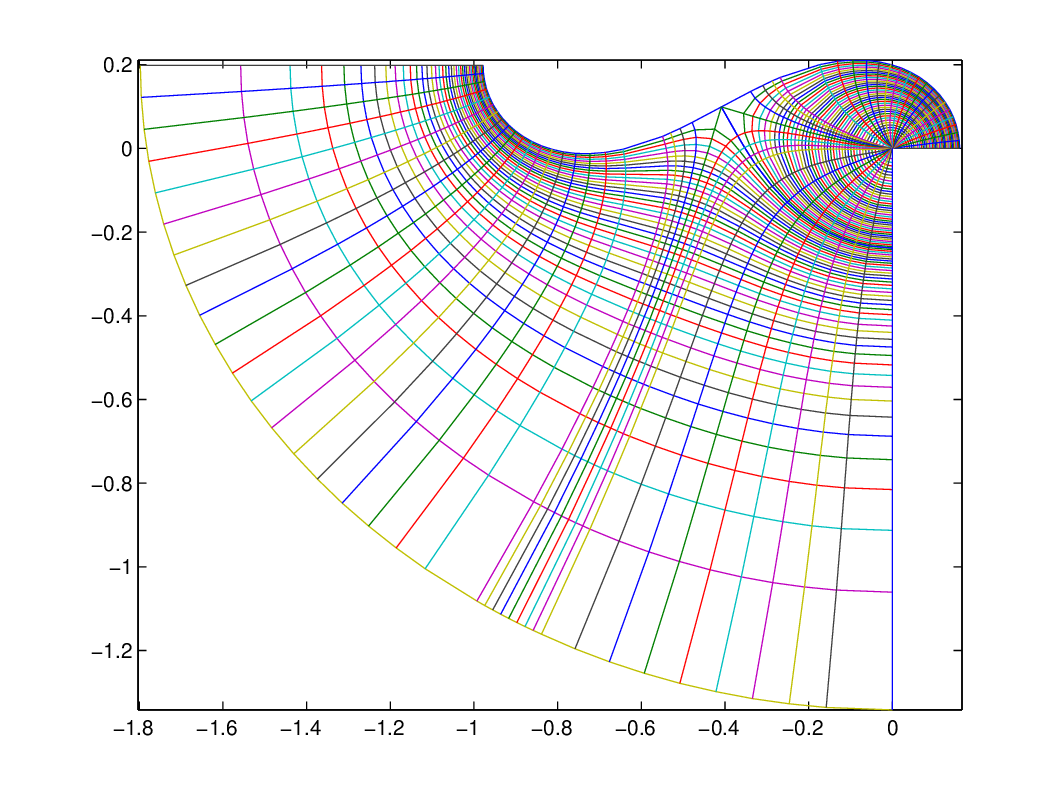}
\caption{Its image under $\zeta$.}
\label{fig:quarterz}
\end{minipage}
\end{figure}

Since $\W-\hz^{-2}\to0$ when $\hz\to0$ then we take $-\int^z\wp=\zeta(z)+\CC/\c$ and find the correct $\CC$ to define $\zeta$. \cref{fig:quarter} shows one quarter of the torus as two mirror copies of \cref{fig:WP_img_cr}, but still avoiding $\infty$ by keeping the original polar coordinates at $\o_1+\o_2$. Its image under $\zeta$ is shown in \cref{fig:quarterz}.

By defining
\BE
   \d:=\lim_{z\to\o_1+\o_2}\biggl(\c\zeta(z)-\frac{1}{z-\o_1-\o_2}\biggl)\label{eq:limit_zeta}
\EE
one gets numerical values of $\d$ for any $\rho\in[0,\pi/2)$. Geometrically speaking, $\d/\c$ is just the (imaginary) level of the horizontal ray that goes to $-\infty$ in \cref{fig:quarterz}. By defining $\eta_1:=\zeta(\o_1)$ and $\eta_2:=\zeta(\o_2)$ it is easy to see that $\d=2|\c|\cdot\Im\,\eta_1=2|\c|\cdot\Im\,\eta_2$. For instance, $\d\cong-0.1132i$ for $\rho=1/2$ because in this case $\c=-0.5697i$ and $\d/\c=0.1987$. By defining $\CC:=-\d+e_3(\o_1+\o_2)$ we finally have
\BE
   \Ze(z)=\CC-e_3z+\c\zeta(z).\label{eq:non_odd}
\EE

Of course, both $\zeta$ and $\Ze$ are odd-functions at $z=0$ and $z=\o_1+\o_2$, respectively. Now we must be careful because one classically defines $-\Et_1$ as $-\Ze$ for the double of the half-period $\hz=-\o_1$ (or $z=\o_2$). Analogously, $-\Et_2$ is $-\Ze$ for the double of the half-period $\hz=-\o_2$ (or $z=\o_1$). Hence, from Equation \cref{eq:non_odd} we have 
\BE
   -\frac{\Et_1}{2}=\CC-e_3\o_2+\c\eta_2=-\d+e_3\o_1+\c\eta_2
   \,\,\,\,{\rm and}\,\,\,\,
   -\frac{\Et_2}{2}=\CC-e_3\o_1+\c\eta_1=-\d+e_3\o_2+\c\eta_1.\label{eq:to_L}
\EE

For $\rho=1/2$ our numerical computations give
\BE
   2(\Et_1\o_2-\Et_2\o_1)=8\c\,\Re\{\eta_1\o_2\}=3.14620i,\label{eq:quasi_pi}
\EE
which is a poor approximation of $\pi i$, namely the exact value of the Legendre relation. This is because in our example we have used a discretization by meshes with $37\times41$ points. For $73\times81$ points we obtain 3.14389 in Equation \cref{eq:quasi_pi} instead, which shows how little improvement one gets with a fourfold refinement.

However, $41\times37$ is our standard mesh dimension, and from Equation \cref{eq:rel_wps} the expression
\BE
   e^{i\rho}-\frac{e_2-e_3}{\c}
\EE
results in $0.0000023362-0.0000042764i$, which is exact up to the 5th decimal place. By using that standard, in \cref{sec:appendix} we shall apply a method which ensures $G$'s accuracy has a maximal error less than $0.002$, and with $143\times321$ points we confirm that our numerical Green's function is reliable until the 4th decimal place for any $\rho\in[-\pi/3,\pi/3]$ with $|T|=1$. Needless to say, numerical analyses cannot go to extreme values, but our methodology remains applicable for a greater range of $\rho$.

Now we shall strongly use \cite{LW} for the computation of the Green's function. Let $T$ be the torus in \cref{fig:scheme}(b) and denote its area by $|T|$. Then~\cite[Eq.(2.6)]{LW} states that
\BE
   8\pi G_{\hz}=\frac{1}{|T|}\int_T\frac{-\W'(\hz)}{\W(\xi)-\W(\hz)}dA,\label{eq:g_odd}  
\EE
and since $\W$ is an even-function, then $G_{\hz}$ is odd. Notice that the concepts of even/odd-functions are indifferently either with respect to the origin or to $\infty$. Moreover, due to Equation \cref{eq:rel_wps} one easily rewrites Equation \cref{eq:g_odd} in terms of the symmetric $\wp$.

\begin{figure}[h]
\center
\includegraphics[scale=0.85]{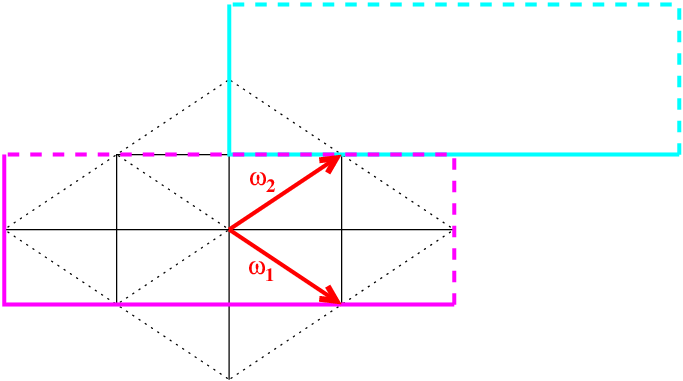}
\caption{The rectangle $\mR$ and the rhombus.}
\label{fig:wp_rec}
\end{figure}

From Equation \cref{eq:g_odd} one sees that $G$ is well-defined in the whole complex plane $\C$. Of course, the sets $2(m\o_1+n\o_2)+T$ will cover $\C$, but from now on we shall take another fundamental domain. With $\a=\Re\,\o_1$ and $\b=\Im\,\o_2$ from Equation \cref{eq:a_and_b}, the following is a reproduction of \cref{fig:scheme}(b) but now we draw an auxiliary rectangle $\mR$ of basis $4\a$ and height $2\b$
centred at the origin. A copy of $\mR$
is also drawn in \cref{fig:wp_rec}, which we call $2\o_2+\mR$. It is obvious that the sets $2(m\o_1+n\o_2)+\mR$ will cover the whole $\C$.

Before going ahead, it is important to notice that many results of \cite{LW} use the following simplification: $\o_1=1$ and $\o_2=\tau=a+ib$. Of course these must be avoided because in our case $\o_2=\ov{\o_1}$. However, we still can profit from some equalities that {\it do not} depend on that simplification. One of them is \cite[Eq.(9.9)]{LW},
\[
   \Ze(\hz)-\Et_1\hz=(\log\vartheta_1(\hz))_{\hz},
\]
which obviously holds if we replace $\hz$ with $z$. For convenience of the reader, here we reproduce the definition of the Jacobi-Riemann theta function
\[
   \vartheta_1(z)=-i\sum_n(-1)^n\bigl(\exp(i\pi\tau)\bigl)^{(n+1/2)^2}\exp{\bigl((2n+1)i\pi z}\bigl).
\]

By the way, let's recall the complex operator
\[
   \frac{\deh}{\deh z}=\frac{1}{2}(\frac{\deh}{\deh x}-i\frac{\deh}{\deh y}).
\]
Now consider~\cite[Eq.(5.1)]{LW}, namely
\BE
   G(\hz)=-\frac{1}{2\pi}\Re\int(\Ze(\hz)-\Et_1\hz)d\hz+\frac{y^2}{2b}+\mathcal{C}(\tau),\label{eq:G_LW}
\EE
where $y=\Im\,\hz=\Im\,z$, $\tau=\o_2/\o_1$, $b=\Im\,\tau$ and $\mathcal{C}$ is a constant determined only by $\tau$. Then it is straightforward to see that Equation \cref{eq:G_LW} implies
\BE
   \frac{\deh G}{\deh z}=\frac{-1}{4\pi}\biggl((\log\vartheta_1)_z+2\pi i\frac{y}{b}\biggl).\label{eq:G_z}
\EE

As explained at the Introduction our numerical $G$ will give $\int_TGdA=-\mathcal{C}(\tau)|T|$ instead of zero. Now we are ready to present our results in the next section.

\section{Experimental results}
\label{sec:experiments}

One sees that Equation \cref{eq:G_z} is precisely the same as in \cite[p.912]{LW}. However, one must be careful because for $\Ze$ in Equation \cref{eq:G_LW} the authors {\it do} use the aforementioned simplification. Anyway, due to Equation \cref{eq:non_odd} for either their approach or ours one has the aforementioned \cref{thm:mainthm}, which we restate here as:

\begin{theorem} Let $\ZZ:\C\to\hat{\R}$ given by
\BE
   \ZZ(\sz)=-\frac{1}{2\pi}\Re\int_0^\sz\c\zeta(z)dz.\label{eq:def_ZZ}
\EE
Then there exist unique constants $A$, $B$, $C$ depending only on $(\a,\b)$ such that $\hat{\ZZ}:\C\to\hat{\R}$ defined by $\hat{\ZZ}(z):=\ZZ(z)+Ax^2+By^2+Cxy$ is a real-analytic doubly-periodic function.
\end{theorem}

\begin{proof} By looking at \cref{fig:wp_rec} one sees three necessary conditions for $\hat{\ZZ}$ to be well-defined along $\ov{\mR}\cap(2\o_2+\ov{\mR})$:
\[
   \ZZ(-2\a-i\b)+A(-2\a)^2+B(-\b)^2+2\a\b C=\ZZ(i\b)+B\b^2,
\]
\[
   0=\ZZ(2\o_2)+A(2\a)^2+B(2\b)^2+4\a\b C,
\]
\[
   \ZZ(-2\a-i\b)+A(-2\a)^2+B(-\b)^2+2\a\b C=\ZZ(2\a-i\b)+A(2\a)^2+B(-\b)^2-2\a\b C.
\]
   
Once again, our symmetric functions will be quite handy. \cref{fig:eagle}(a) outlines the fundamental domain of $\wp$, depicted in \cref{fig:quarter}, together with subsequent reflection in the segment $\ov{MJ}$, $180^\circ$-rotation around $N$ and reflection in $\ov{BM}$. As a matter of fact, these points were named after the eagle-profile in \cref{fig:quarterz}: $J$ugular notch, $M$outh, $B$eak, $N$ape, $S$houlder and $W$ing. \cref{fig:eagle}(b) shows what happens to the eagle-profile, which ends up in a translation of $2\eta_2$. Since $\ZZ$ is even, $\zeta$ odd and $\zeta(2\o_2)=2\eta_2$, with the change of variables $z\mapsto z+2\o_2$ it is not difficult to obtain $\ZZ(2\o_2)=2\ZZ(\o_2)-\frac{1}{2\pi}\Re\{\c\eta_2\o_2\}$.

\begin{figure}[h!]
\center
\includegraphics[scale=0.7]{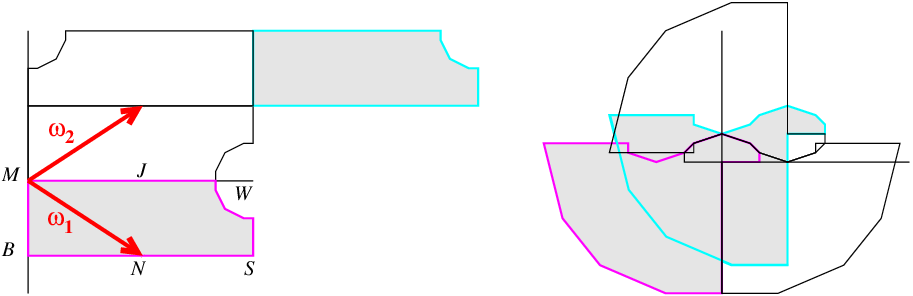}
\centerline{(a)\hspace{6.5cm}(b)}
\caption{(a) Reflections in the fundamental domain, (b) translation by $2\eta_2$.}
\label{fig:eagle}
\end{figure}
  
Now a simple computation gives
\BE
   A=\frac{\ZZ(i\b)-\ZZ(2\a+i\b)}{4\a^2},
   B=\frac{\ZZ(2\a+i\b)-\ZZ(i\b)-2\ZZ(\o_2)+0.5\Re\{\c\eta_2\o_2\}/\pi}{4\b^2},C=0.
   \label{eq:three_consts}
\EE

By the way, the change of variables $z\mapsto z+2\o_2$ also leads to the general formula
\BE
   \ZZ(z+2\o_2)=\ZZ(z)-\frac{1}{\pi}\Re\{\c\eta_2z\}+
   2\ZZ(\o_2)-\frac{1}{2\pi}\Re\{\c\eta_2\o_2\},\label{eq:gen_ZZ}
\EE
which we apply for $z=-\o_2$ and $z=-2\a-i\b$ to get
\[
   \ZZ(\o_2)=-\frac{1}{4\pi}\Re\{\c\eta_2\o_2\}
   \,\,\,\,\,\,{\rm and}\,\,\,\,\,\,
   \ZZ(i\b)-\ZZ(2\a+i\b)=\frac{\a}{\pi}\Re\{\c\eta_2\},
\]
respectively. Now we recall that $\c=-i|\c|$ in order to simplify Equations \cref{eq:three_consts} and \cref{eq:gen_ZZ} as
\BE
   A=\frac{|\c|}{4\pi\a}\Im\eta_2,\,\,\, 
   B=\frac{|\c|}{4\pi\b}\Re\eta_2,\,\,\,C=0,\,\,\,
   \ZZ(z+2\o_2)=\ZZ(z)-\frac{1}{\pi}\Re\{\c\eta_2z\}-
   \frac{1}{\pi}\Re\{\c\eta_2\o_2\},\label{eq:simpl_formls}
\EE
respectively. At this point an easy computation gives
\[
   \hat{\ZZ}(z+2\o_2)\stackrel{\rm def}{=}
   \ZZ(z+2\o_2)+A(x+2\a)^2+B(y+2b)^2\stackrel{\cref{eq:simpl_formls}}{=}
\]
\[
   \ZZ(z)+Ax^2+By^2+\frac{|\c|}{\pi}\Re\{i\eta_2z\}+4\a Ax+4\b By+
   \frac{|c|}{\pi}\Re\{i\eta_2\o_2\}+4\a^2A+4b^2B
   \stackrel{\cref{eq:simpl_formls}}{=}
\]
\[
   \ZZ(z)+Ax^2+By^2=\hat{\ZZ}(z),
\]
and similarly $\hat{\ZZ}(z+2\o_1)=\hat{\ZZ}(z)$. Hence $\hat{\ZZ}(\cdot+2m\o_1+2n\o_2)=\hat{\ZZ}\,\forall m,n$.
\end{proof}

It is important to notice \cref{cor:a}, which we restate here as:

\begin{corollary}
$A$ and $B$ from Equation \cref{eq:simpl_formls} result in $4(A+B)=1/|T|$.
\end{corollary}

\begin{proof} We recall that $\d=2|\c|\cdot\Im\eta_1=2|\c|\cdot\Im\eta_2$, hence Equation \cref{eq:to_L} implies
\[
   2(\Et_1\o_2-\Et_2\o_1)=4\d(\o_2-\o_1)+4\c(\eta_1\o_1-\eta_2\o_2)=
   8|\c|\Im\eta_1\cdot 2\b i-4|\c|(\eta_1\o_1-\eta_2\o_2)i.
\]

But $\eta_1\o_1-\eta_2\o_2=2\a\Re\eta_1+2\b\Im\eta_1$, whence $2(\Et_1\o_2-\Et_2\o_1)=4|\c|(\eta_2\o_1-\eta_1\o_2)i$, and by Legendre's relation we have
\BE
   4|\c|(\eta_2\o_1-\eta_1\o_2)=\pi.\label{eq:my_L}
\EE

Since $\eta_2=-\ov{\eta}_1$, Equations \cref{eq:simpl_formls} and \cref{eq:my_L} give
\[
   4(A+B)=\frac{|\c|}{\pi\a\b}(\b\Im\eta_1-\a\Re\eta_1)=
   \frac{|\c|}{\pi\a\b}\cdot\frac{\pi}{8|\c|}=
   \frac{1}{8\a\b}=\frac{1}{|T|}.
\]

\end{proof}

\begin{figure}[h!]
\center
\includegraphics[scale=0.3]{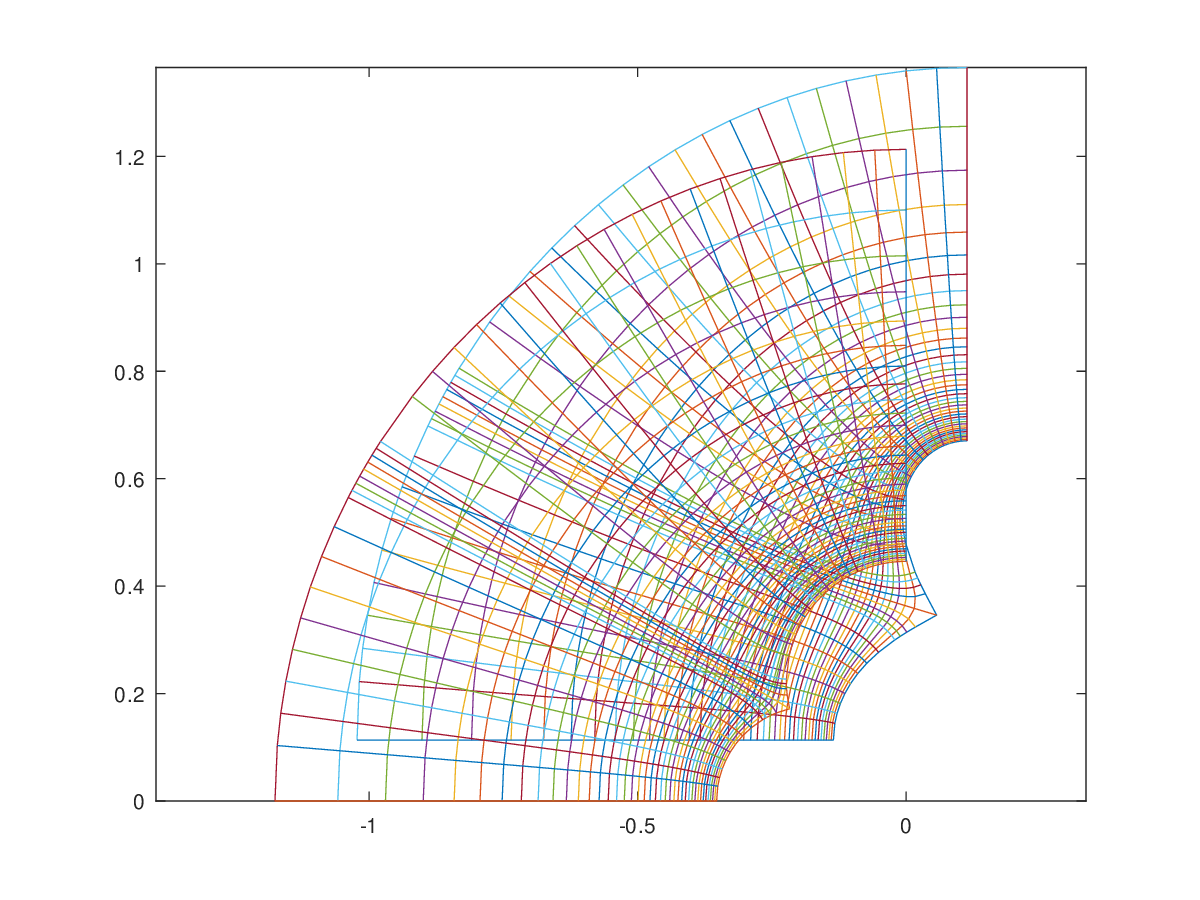}
\caption{Overlap between the numerical functions $\c\zeta-\d$ and $0.5/(z-2\a)$.}
\label{fig:numerical_zeta}
\end{figure}

At this point we ought to comment that our theoretical results are closely accompanied by helpful numerical simulations, but these can be tricky. Indeed, \cref{fig:WP_img_cr} was taken as a subdomain of $\wp$, which maps it to the right half of the unit disk in \cref{fig:fund_rectg_cr}. But we could have done it with any dilation of \cref{fig:WP_img_cr} as well, which will then affect the {\it numerical} functions. Indeed, from Equation \cref{eq:limit_zeta} we see that $\c\zeta-\d$ must approach $1/(z-2\a)$, and \cref{fig:numerical_zeta} shows a comparison between the former and {\it half} the latter, both restricted to the right half of \cref{fig:quarter}.
  
In point of fact, our numerical functions are already the theoretical ones divided by two. Namely, from Equation~\cref{eq:classical_WP} we have been working with $\W(z/\sqrt{2})/2$. In order to confirm it, our program {\tt rh} mentioned in \cref{sec:alg} prints the maximal fitting gap between the spatial curves $BN$ and $SN$ (this one rotated $180^\circ$ around $N$), which is $1.6698\cdot10^{-9}$ for $\rho=1/2$. Compare \cref{fig:eagle}(a) with \cref{fig:graph} for a visualisation of these curves. Therefore, in \cref{thm:def_G} we have decided to write $G$ with the factor 2 in order to make it compatible with our numerical simulations.

Now we are ready to establish our Green function presented in \cref{thm:def_G}, which we restate here as

\begin{theorem}
Consider $\ZZ$ as in \cref{thm:mainthm} and $G=2\ZZ+Ax^2+By^2$. Then $G$ is the Green's function with $G(0)=0$ and $\Delta G=-\sum_{m,n}\delta(2\a+2m\o_1+2n\o_2)+1/|T|$ in $\C$.
\end{theorem}

\begin{proof} In view of what we have just explained, together with the fact that ours is the non-negative Green function with $G(0)=0$, Equation \cref{eq:G_LW} is rewritten as:
\BE
   G(\hz)=-\frac{1}{2\pi}\Re\int(2\Ze(\hz)-\Et_1\hz)d\hz+\frac{y^2}{2b}.\label{eq:new_G_LW}
\EE
Now we use \cref{eq:non_odd} and \cref{eq:new_G_LW} to get:
\[
   G(\sz)=-\frac{1}{\pi}\Re\int_0^\sz\c\zeta(z)dz+p(x,y),
\]
where $p(x,y)$ is a quadratic polynomial and $p(0,0)=0$. But $G$ is an even function, whence $p_x(0,0)=p_y(0,0)=0$. Therefore, from \cref{thm:mainthm} we have $G=2\ZZ+Ax^2+By^2$, and since in~\cite{LW} the authors prove that Equation \cref{eq:new_G_LW} is the Green's function, then $\Delta G=-\sum_{m,n}\delta(2\a+2m\o_1+2n\o_2)+1/|T|$ in $\C$. In the next section we give numerical evidence that $G$ is non-negative.
\end{proof}
 
\section{Conclusions}
\label{sec:conclusions}

We have applied symmetric elliptic functions to obtain the numerical non-negative Green's function on any flat rhombic torus $T$. Our work focus on the methodology, which circumvents truncation errors by forgoing infinite series. The other numerical errors are controlled through changes of variables and refinement comparison, this one discussed in \cref{sec:appendix}. The method guarantees an accuracy of four decimal digits with $|T|=1$, which always takes less than 5s to be confirmed with Matlab R2016a on our architecture. We also avoid exceedingly many arithmetic operations, because these will accumulate and propagate machine error.

\appendix
\section{Accuracy of the numerical $G$}
\label{sec:appendix}

As commented right after Equation \cref{eq:quasi_pi}, we can adopt the strategy of quadrupling our standard mesh and then check which decimal places of the numerical $G$ will remain invariant after the refinement. Of course, this must be done for many values of $\rho$, but we can restrict our analysis to what is depicted in \cref{fig:quarter} and \cref{fig:quarterz}. Actually, our starting domain is illustrated in \cref{fig:WP_img}(a), which we have full control of. This one will be refined in order to obtain the numerical $\GG$, whose decimal places will be compared with $G$'s. 

Since we are going to fix $|T|=1$, both $\GG$ and $G$ values will be previously multiplied by $k^2$, where $k:=2\sqrt{2\a\b}$. Indeed, from Equation \cref{eq:def_G} we see that $\W$ should then be defined for $\o_1/k$ and $\o_2/k$, hence a re-scaling of $\C$ by $1/k$. Now, for any positive real $\kappa$ we consider the one-dimensional probability distribution $\delta_\kappa:\R\to\R$, given by $\delta_\kappa(t)=\kappa\exp(-(\kappa t)^2)/\sqrt{\pi}$. Then the change of variables $t\mapsto t/k$ gives
\[
   \int_{-\infty}^\infty\frac{\kappa\exp(-(\kappa t/k)^2) }{k\sqrt{\pi}}dt=1
\]
for the distribution $\delta_{\kappa/k}$. Since $\delta(x,y)=\delta(x)\cdot\delta(y)$, then
\[
   \Delta(k^2G) = -\delta+1\hspace{1cm}{\rm with}\hspace{1cm}\int_T(k^2G)dA=k^2\int_TGdA=0.
\]
By the way, this is the reason why $\mathcal{C}$ in Equation \cref{eq:G_LW} only depends on $\tau$. The same conclusion can also be reached from Equation \cref{eq:g_odd}, and alternatively from the fact that $\ln|z/k|=\ln|z|-\ln|k|$.

After downloading {\tt rh.zip} from our website, extract it in a folder and at the Matlab prompt enter {\tt rh} (for the surface) or {\tt rh4} to check the accuracy. Either will start as follows:
\bigskip
\begin{verbatim}
For numerics we must avoid extremes
Thus give rho in [-pi/3,pi/3]:
\end{verbatim}
\bigskip

Although {\tt rh4} uses a finer mesh than {\tt rh}, this one also includes the option to check $G$'s accuracy. In all of our tests with {\tt rh} the worst accuracy fails only at the 3rd decimal place, but $0.001813$ was the maximal difference between $\GG$ and $G$ (attained for $\rho=0$). However, {\tt rh4} does {\it not} draw the surface and it confirms that the accuracy always reaches at least four decimal places.

Some information is printed by {\tt rh}, like $\eta_1$, $\eta_2$ and the maximal fitting gap mentioned in \cref{sec:experiments}. By choosing to check precision {\tt rh} invokes {\tt rheavy}, which confirms that $\GG\ge0$ and gives the accuracy by printing the maximal absolute error, namely $\max|\GG-G|$. For instance, it is $0.001467$ for $\rho=1/2$.

This value drops to $0.000043$ if we run {\tt rh4}, because it works with a $143\times321$-mesh for $G$, as explained in \cref{sec:alg}. Here it is compared with a refinement $\GG$ of $143\times641$ points. Another two meshes of $37\times41$ and $73\times81$ points were also cited in \cref{sec:alg}, but as one sees we do not always double both dimensions in our programs. Actually the first dimension also depends on $\rho$, because the discretization must be denser at $e^{i\rho}$, as depicted in \cref{fig:WP_img}(a). This point causes a singularity that was removed with Equation \cref{eq:get_rid}, but still requires refinement to render a smooth surface. That is why a trace appears towards each pinnacle in \cref{fig:graph}.

In order to estimate the constant $\int_TGdA$ we cannot remove singularities by changing variables, as done in \cref{sec:alg}. There we had algebraic formulae, but now our $G$ is numerical. Hence we must integrate $G$ over the rectangle in \cref{fig:quarter}, and at $z=2\a$ compute the logarithmic integral separately. Now
\[
   -\int\!\!\!\!\int_{B_\eps(0)}\ln|\hz|dxdy=\frac{\pi\eps^2}{2}(1-2\ln\eps)
\]
for a circle of radius $\eps$ around $z=2\a$. We can take the average of the $G$-values along the corresponding circumference, and since that integral will be dominant, $\int_TGdA$ is calculated by adding up the three parts. However, we decided not to implement this strategy because it penalizes the accuracy without a clear error control.

\section*{Acknowledgements}
We would like to thank Prof Frank Baginski at The George Washington University for his helpful assistance in the elaboration of this manuscript.

\bibliographystyle{siamplain}
\bibliography{axel}
\end{document}